\documentclass{amsart}
\usepackage{mathrsfs}

\usepackage{amsmath, amsfonts, amssymb,amsthm}
\usepackage[mathscr]{eucal}

\usepackage{color}
\usepackage{graphicx}
\usepackage[T1]{fontenc}
\usepackage[sc]{mathpazo}
\usepackage[all]{xy}
\usepackage{hyperref}
\usepackage{chngpage}
\usepackage{array}

\usepackage{enumerate}

\newtheorem{theorem}{Theorem}[section]

\newtheorem{lemma}[theorem]{Lemma}

\definecolor{alert}{rgb}{0.8,0,0}

\title[volume of hypersurface]
{On the volume of locally conformally flat 4 dimensional closed hypersurface}
\thanks{This work is partially supported by National Nature Science Foundation
 of China (Grant No. 11601442) and
Fundamental Research Funds for the Central Universities (Grant No. 2682016CX114, WK0010000055).}

\author{Qing Cui}
\address{School of Mathematics,
Southwest Jiaotong University, 611756
Chengdu, Sichuan, China}
\email{qingcui@impa.br}

\author{Linlin Sun}
\address{School of Mathematics Sciences,
University of Science and Technology of China,
230026 Hefei, Anhui, China}
\email{sunll@ustc.edu.cn}



\begin{document}
\maketitle

\begin{abstract}
Let $M$ be a 5 dimensional Riemannian manifold with $Sec_M\in[0,1]$,
$\Sigma$ be a locally conformally flat closed hypersurface in $M$ with
mean curvature function $H$. We prove that, there exists
$\varepsilon_0>0$, such that
\begin{align}\label{V1}
\int_\Sigma (1+H^2)^2 \ge \frac{4\pi^2}{3}\chi(\Sigma),
\end{align}
provided $\vert H\vert \le \varepsilon_0$, where $\chi(\Sigma)$ is the
Euler number of $\Sigma$. In particular,
if $\Sigma$ is a locally conformally flat minimal hypersphere in
$M$, then $Vol(\Sigma) \ge 8\pi^2/3$, which partially answer
a question proposed by Mazet and Rosenberg \cite{Ma&Rosen}.
%For
%an $(n+1)-$ dimensional rotationally symmetric Riemannian manifold
% $M$, we show that an immersed hypersurface $\Sigma$ is locally
% conformally
% flat if and only if ($n-1$) of the principal curvatures of $\Sigma$
%are the same, which is a generalization of Cartan's result
% \cite{Cartan}.
Moreover, we show that if $M$ is (some special but large class) rotationally symmetric,
%5-manifold with $Sec_M\in [0,1]$, and $\Sigma$ is a locally conformally
%flat closed hypersurface
% with mean curvature $H$,
 the
inequality (\ref{V1})
holds for all $H$.
\end{abstract}
\maketitle
% ----------------------------------------------------------------
\section{Introduction}

Let $M$ be a  2-sphere with a smooth Riemannian metric
such that the curvature is between 0 and 1. It is known (see
\cite{Kligen} or \cite{P}) that the length of an embedded closed geodesic in
$M$ is at least 2$\pi$, which is the length of the
standard circle in Euclidean plane. When $M$ is a
Riemannian 3-manifold with sectional curvature between 0 and 1, one
can easily apply Gauss equation and Gauss-Bonnet theorem to obtain
that an embedded minimal sphere $\Sigma$ in $M$ has area at least
$4\pi$, that is
$$
4\pi = \int_{\Sigma} Sec_\Sigma =\int_{\Sigma} R_{1212}
= \int_{\Sigma} \left(\overline{R}_{1212} - \frac{1}{2}\vert A\vert^2\right)
\le \int_\Sigma \overline{R}_{1212} \le Vol\left(\Sigma\right),
$$
where $R$ and $\overline{R}$ denote the curvature tensors of $\Sigma$
and $M$, $A$ denotes the second fundamental form of $\Sigma$ in $M$.

In
\cite{Ma&Rosen}, Mazet and Rosenberg study the equality case and get a
rigidity theorem for $M$.
The authors also put forward two
very interesting questions, one of them is, if
$M$ is an ($n$+1)-Riemannian manifold with $Sec_M\in[0,1]$,
does an embedded minimal hypersphere
(i.e., minimal hypersurface diffeomorphic to the
standard Euclidean $n$-sphere $\mathbb S^n$) has
volume at least the volume of $\mathbb S^n$?
In 1974, Hoffman and Spruck \cite{HS} studied the isoperimetric
inequality and showed, if $M$ is a simply connected
Riemannian ($n+1$)-manifold with $Sec_M\in [1/4, 1]$, then any closed minimal
hypersurface has at least the volume of $\mathbb S^n$.
Therefore, if the answer of the Mazet and Rosenberg's question
is true, it can be seen as a generalization (with topological restricted) of Hoffman
and Spruck's result. We would like to point out that, if
$Sec_M\in[0,1]$, the topological restriction on $\Sigma$ is necessary. Actually, given
$\varepsilon>0$, let
 $\Sigma$ be a flat $n-$torus  with $Vol(\Sigma) \le \varepsilon$ (which can be done
 by passing a dilation), then $\Sigma$ is a totally geodesic closed hypersurface embedded
 in
 $\Sigma\times\mathbb R$ whose sectional curvature is 0.

Note that in the
case of $n=2$, every surface admits an isothermal coordinates and
therefore is locally conformally flat. It is seems natural to add
 the condition "locally
conformally flat" on the hypersurface in high dimensional case.
In this paper, we focus our attention
on the case of $n=4$, pose the assumption that $\Sigma$ is locally conformally
flat, and partially answer the question proposed by
Mazet and Rosenberg. Actually, we get a more general result as follows.
\begin{theorem}\label{1}
Let $M$ be a 5-dimensional Riemannian manifold with
$Sec_M\in[0,1]$, and $\Sigma$ be an embedded
locally conformally flat closed hypersurface in
$M$ with mean curvature function $H$.
Then we have
\begin{align}\label{vol}
\int_{\Sigma} \left( (1+H^2)^2 +\vert H\vert f(\vert H\vert)\right)\ge  \frac{4\pi^2}{3}\chi(\Sigma) \ ,
\end{align}
where $f$ is a nonnegative function defined in Section 2, and $\chi(\Sigma)$ is the
Euler number of $\Sigma$.

Moreover, there exists $\varepsilon_0>0$, such that if
$\vert H\vert \le \varepsilon_0$, we obtain,
\begin{align}\label{vol1}
\int_{\Sigma} \left(1+H^2\right)^2\ge  \frac{4\pi^2}{3}\chi(\Sigma).
\end{align}
The equality holds if and only if
the mean curvature $H$ is constant, $\Sigma$ is totally umbilic
and isometric to $\mathbb S^4\left(\frac{1}{1+H^2}\right)$.
\end{theorem}

As an immediate
corollary of theorem \ref{1}, the following result partially
answer the question proposed by Mazet and Rosenberg.

\begin{theorem}
Let $M$ be a 5-dimensional Riemannian manifold with
$Sec_M\in[0,1]$, and $\Sigma$ be an embedded
locally conformally flat  minimal hypersphere in
$M$.
Then
\begin{align*}
Vol(\Sigma)\ge  \frac{8\pi^2}{3} =Vol\left(\mathbb S^4\right).
\end{align*}
The equality holds if and only if
$\Sigma$ is totally geodesic
and isometric to $\mathbb S^4$.
\end{theorem}
This paper is organized as follows. In Section 2, we list
 some notations and known formulas, and give the proof
 of Theorem \ref{1}. In Section 3,
% we prove a generalization
% of Cartan's result about the relation of locally conformally
% flatness and the principal curvatures, see Lemma \ref{principal}.
% Using this lemma,
 we deal with a special case when $M$ is rotationally
 symmetric and get the lower bound volume for all $H$, see Theorem \ref{rotsym}.

\section{Preliminary and proof of Theorem \ref{1}}

Let $(M,\bar{g})$ be an $n+1$ dimensional Riemmannian manifold,
and $(\Sigma,g)$ be a hypersurface isometric immersed in $M$.
If there is no ambiguity,
$\langle \cdot \ ,\ \cdot\rangle$ will denote both $\bar{g}$ and $g$.
Let $\overline{\nabla}$ and $\nabla$ be the Levi-Civita connection
induced by metric $\bar{g}$ and $g$ respectively. Let $R$ be the curvature
tensor on $\Sigma$ defined by, for all $X,Y,Z,W\in\mathfrak{X}(T\Sigma)$,
$$
R(X,Y,Z,W) = \langle R(X,Y)Z,W\rangle,
$$
where $R(X,Y) =-\nabla_X\nabla_Y +\nabla_Y\nabla_X +\nabla_{[X,Y]}$. Also let
 $\overline{R}$ be the curvature tensor on $M$ which is defined similarly.

Let
$e_1, \cdots, e_n$ be a local orthonormal frame on $\Sigma$. For all
$1\le i, j, k, l\le n$, write
$$
R_{ijkl} = R(e_i, e_j, e_k, e_l),\quad
\overline{R}_{ijkl} = \overline{R}(e_i, e_j, e_k, e_l).
$$
The sectional curvature will be
$$
Sec_M(e_i\wedge e_j) = R_{ijij},\quad
Sec_{\overline{M}}(e_i\wedge e_j) = \overline{R}_{ijij}.
$$
Let $A$ be the second fundamental form of $\Sigma$
in $M$, $h_{ij} = \langle A(e_i), e_j\rangle$ be the coefficients of
$A$.
Then the Gauss equation can be written as
\begin{align*}
R_{ijkl} = \overline{R}_{ijkl} +h_{ik}h_{jl} - h_{il}h_{jk}.
\end{align*}

We also denote by $\mathbb S^n $ be the standard unit $n$ sphere in $n+1$ Euclidean
space and by $\mathbb S^n(r)$ be the round $n$-sphere with radius $r$.
Now we will prove Theorem \ref{1}.

\noindent
{\it Proof of Theorem \ref{1}.} The Gauss-Bonnet-Chern
formula for a closed 4-manifold $\Sigma$ is
(see \cite{Avez} or \cite{Besse})
\begin{align}\label{GBC1}
4\pi^2\chi \left(\Sigma\right) = \int_{\Sigma}\left( \frac{S^2}{12}
-\frac{\vert Ric\vert^2}{4} + \frac{\vert W\vert^2}{8}\right),
\end{align}
where $\chi(\Sigma)$ is the Euler characteristic of $\Sigma$,
$S$ is the scalar curvature, $Ric$ is the Ricci tensor and
$W$ is the Weyl tensor.  It is well known that,
 when dimension grater than 3,
locally conformally flatness equivalent to Weyl tensor vanishing.
Therefore, to prove the first part of the
 theorem, it is sufficient to prove,
pointwisely,
\begin{align*}
\text{{\bf Claim.}}\quad\quad\quad\quad
Q:= \frac{S^2}{12} - \frac{\vert Ric\vert^2}{4} \le 3(1+H^2)^2
+ 3 \vert H\vert f(\vert H\vert).\quad\quad\quad\quad\quad\quad
\end{align*}
Next we will consider our problem at one point
$p\in \Sigma$ (in the calculations, we omit the
 letter "$p$" for
simplicity). Throughout this proof, $i, j, k, l$
will range from 1 to 4 if there is no special explanation.

Firstly, by the Gauss equation, we obtain,
\begin{align}\label{scalar}
S^2 &= \left( \sum_{i,j} R_{ijij}\right)^2
=\left( \sum_{i,j}\overline{R}_{ijij} +16H^2-\vert A\vert^2\right)^2
=\left(\sigma+12H^2-\vert\mathring{A}\vert^2\right)^2\\
&= \sigma^2 + 144H^4 +\vert\mathring{A}\vert^4
+24\sigma H^2 -2\sigma\vert\mathring{A}\vert^2-24H^2\vert\mathring{A}\vert^2,\notag
\end{align}
where
$$\sigma:=\sum_{i,j}\overline{R}_{ijij}, \quad \text{and}\quad
 \mathring{A} := A - HI,$$
 i.e., $\mathring{A}$ is the traceless
part of $A$.

For simplicity, let $e_1, e_2, e_3, e_4$ be the principal
directions at the point $p$, and $\lambda_1, \lambda_2,
\lambda_3, \lambda_4$ be the corresponding principal curvatures,
we have
\begin{align*}
\vert Ric\vert^2 = \sum_{i,j} \left(\sum_{k} R_{ikjk}\right)^2
= \sum_{i,j} \left( \sum_{k}\overline{R}_{ikjk} +4\delta_{ij} \lambda_i H
-\delta_{ij}\lambda_i \lambda_j \right)^2
%&= \sum_{i,j} \left(\sum_{k} \overline{R}_{ikjk}\right)^2
%- 2\sum_{i}\left(\lambda_i^2 \sum_{k} \overline{R}_{ikik}\right)
%+\sum_{i} \lambda_i^4\notag
\end{align*}
For simplifying $\vert Ric\vert^2$, we need to introduce some notations as follows:
\begin{align*}
a_{ij} := \sum_{k} \overline{R}_{ikjk}, \quad
\mathring{a}_{ij} := \sum_{k} \overline{R}_{ikjk} - \frac{\sigma}{4}\delta_{ij}.
\end{align*}
Note that $\sigma$ is the trace of $(a_{ij})$
and $(\mathring{a}_{ij})$ is the traceless part of $(a_{ij})$.
Using these notations, we get
\begin{align}\label{Ric1}
\vert Ric\vert^2=&\sum_{i,j}\left(a_{ij}+4\delta_{ij}\lambda_i H
-\delta_{ij}\lambda_i\lambda_j\right)^2 \\
=& \sum_{i,j} a_{ij}^2 +16H^2\vert A\vert^2 +\sum_{i} \lambda_i^4\notag\\
&+8H \sum_i \lambda_i a_{ii} - 2\sum_{i}\left(\lambda_i^2 a_{ii}
\right) -8H\sum_i\lambda_i^3
\notag\\
=& \frac{\sigma^2}{4} + \vert\mathring{a}\vert^2
+16H^2\vert \mathring{A}\vert^2 +64H^4 +\sum_{i} \lambda_i^4\notag\\
& - 2\sum_{i}\left(\lambda_i^2 -4H\lambda_i\right) a_{ii}
 -8H\sum_i\lambda_i^3,
\notag
\end{align}
where $\vert \mathring{a}\vert^2 = \sum_{i,j} \mathring{a}_{ij}^2$.
Next we set $\mu_i = \lambda_i-H$ which are the eigenvalue of $\mathring{A}$.
Then by a direct computation, we have
\begin{align}\label{mu4}
\sum_i \lambda_i^4 &= \sum_i\mu_i^4 +12H^4 -6H^2\vert A\vert^2 +4H\sum_i\lambda_i^3\\
&= \sum_i\mu_i^4 -12H^4 -6H^2\vert \mathring{A}\vert^2 +4H\sum_i\lambda_i^3,\notag
\end{align}
and
\begin{align}\label{mu3}
\sum_i \lambda_i^3 &= \sum_i \mu_i^3 +4H^3 +3H\vert\mathring{A}\vert^2.
\end{align}

Combined (\ref{scalar}), (\ref{Ric1}), (\ref{mu4}) and (\ref{mu3})
we obtain
\begin{align}\label{SR}
Q
 =\frac{1}{12}&\left(\frac{\sigma^2}{4} + 6\sigma H^2 +36H^4
 +\vert \mathring{A}\vert^4 - 3\sum_i\mu_i^4
 +6\sum_i \mu_i^2\left(a_{ii}-\frac{\sigma}{3}\right)\right.\\
 &\left.\ \ -12H\sum_i \mu_i a_{ii} -18H^2\vert \mathring{A}\vert^2
 +12H\sum_i \mu_i^3 -3 \vert\mathring{a}\vert^2\right).\notag
\end{align}

We will divide our proof of the Claim into two main cases. We will see that the locally conformally flatness
will play a key role in the estimate of $Q$.\\\\
{\bf (i)} At the point $p$, $\vert \mathring{A}\vert^2 (p) \le 12 + 24 H^2(p)$.

By the Gauss equation,
$$
S = \sigma + 12H^2 -\vert \mathring{A}\vert^2,
$$
where $\sigma$ is defined in (\ref{scalar}). Since $0\le \sigma\le 12$,
on one hand
$$
S = \sigma + 12H^2 -\vert \mathring{A}\vert^2 \le 12 + 12H^2.
$$
On the other hand,
$$
S = \sigma + 12H^2 -\vert \mathring{A}\vert^2\ge 12H^2
-(12 + 24H^2) = -12(H^2+1).
$$
The above two inequalities yield $S^2 \le 144(1+H^2)^2$.
Therefore,
\begin{align*}
Q = \frac{S^2}{12}-\frac{\vert Ric\vert^2}{4}
= \frac{S^2}{48}-\frac{\vert E\vert^2}{4}
\le \frac{S^2}{48} \le  3(H^2+1)^2,
\end{align*}
where $E$ is the traceless
part of the Ricci tensor, namely the Einstein tensor.\\\\
{\bf (ii)} At the point $p$,
$\vert \mathring{A}\vert^2 (p) \ge 12 + 24 H^2(p)$.

The proof of this case is more difficult than case (i). To prove
the claim, we need to estimate $Q$ by using the equality
(\ref{SR}).

First note that for a fixed $i$, the term $ a_{ii}
-\frac{\sigma}{3}$ is bounded above by $1$. We take $i=1$ for example:
\begin{align}\label{R}
& a_{11} -\frac{\sigma}{3} = \sum_k\overline{R}_{1k1k}
-\frac{\sigma}{3} =\sum_k\overline{R}_{1k1k}
-\frac{1}{3}\sum_{i,j} \overline{R}_{ijij}\\
&=\frac{1}{3}\left(\overline{R}_{1212} + \overline{R}_{1313}
+ \overline{R}_{1414}\right)
-\frac{2}{3}\left(\overline{R}_{2323} + \overline{R}_{2424}
+ \overline{R}_{3434}\right)\le 1\notag
\end{align}
where we have used the curvature condition that $0\le
\overline{R}_{ijij}\le 1$ for all $i\neq j$.

By a direct computation, we have
\begin{align}\label{lambda}
\sum_i \mu_i^4 &= \sum_{i=1}^3 \mu_i^4 + \left(\sum_{i=1}^3 \mu_i\right)^4\\
&=\frac{1}{2} \left(\sum_{i=1}^3 \mu_i^2 +\left(\sum_{i=1}^3 \mu_i\right)^2 \right)^2
+ 4\mu_1\mu_2\mu_3\left(\sum_{i=1}^3 \mu_i\right)\notag\\
&=\frac{1}{2}\vert \mathring{A}\vert^4 - 4\prod_i \mu_i
:=\frac{1}{2}\vert \mathring{A}\vert^4 - 4\mathcal{K},\notag
\end{align}
where $\mathcal{K} = \prod_i\mu_i$ is the Gauss-Kronecker curvature
of $\mathring{A}$.

Obeserve that $\sum_i \mu_i =0$, we get
\begin{align}\label{HA}
&-12H\sum_i \mu_i a_{ii}  -3 \vert\mathring{a}\vert^2\\
&=-12H\sum_i \mu_i \mathring{a}_{ii}  -3 \vert\mathring{a}\vert^2\notag\\
&=-3\sum_i \left( \mathring{a}_{ii}^2 -4H \mu_i \mathring{a}_{ii} +4H^2\mu_i^2\right)
+12H^2\vert\mathring{A}\vert^2 -3\sum_{i\neq j} \mathring{a}_{ij}^2\notag\\
&\le 12H^2\vert \mathring{A}\vert^2.\notag
\end{align}

Combined (\ref{SR}), (\ref{R}), (\ref{lambda}),  (\ref{HA}) and the fact that
$0\le\sigma\le 12$, we have
\begin{align}\label{Q1}
Q\le &3\left(1+H^2\right)^2\\
 &+\frac{1}{12}\left(-\frac{1}{2} \vert \mathring{A}\vert^4
+12H\sum_i \mu_i^3 + 6(1-H^2)\vert \mathring{A}\vert^2
 +12\mathcal{K} \right).\notag
\end{align}

Next we will take the Weyl tensor into consideration.
 The Weyl tensor defined in a coordinate
chart is given by (see e.g. \cite{Aubin}, p117)
\begin{align}
W_{ijkl} = &R_{ijkl} - \frac{1}{2}
\left(R_{ik}g_{jl} - R_{il} g_{jk}
+ R_{jl} g_{ik} - R_{jk} g_{il}\right)\\
&+ \frac{S}{6} \left(g_{jl} g_{ik} - g_{jk} g_{il}\right),\notag
\end{align}
where $R_{ij} = \sum_k R_{ikjk}$ is the Ricci tensor.
Therefore, when $i\neq j$,
we have,
\begin{align*}
W_{ijij} = &R_{ijij} - \frac{1}{2}
\left(R_{ii}
+ R_{jj}\right)+ \frac{S}{6}.
\end{align*}
Now we fix $i=1$ and $j=2$ for example, and get
\begin{align*}
&\frac{S}{6}-W_{1212}\\
&= \frac{1}{2} \left( R_{11} + R_{22}\right) - R_{1212}\\
&= \frac{1}{2}\left(R_{1313} + R_{1414} + R_{2323} +R_{2424}\right)\\
&= \frac{1}{2}\left(R_{1212}+R_{1313} + R_{1414} + R_{2323} +R_{2424}
+R_{3434}\right) - \frac{1}{2}\left(R_{1212}
+R_{3434}\right)\\
&= \frac{S}{4} - \frac{1}{2}\left(R_{1212}
+R_{3434}\right).
\end{align*}
As a consequence, for $\{i,j,k,l\} = \{1,2,3,4\}$, we obtain,
\begin{align}\label{S}
\frac{S}{6} &= R_{ijij} +R_{klkl}-2W_{ijij}\\
=&\overline{R}_{ijij} +\overline{R}_{klkl} -2W_{ijij}
+\lambda_i\lambda_j +\lambda_k\lambda_l\notag\\
=&\overline{R}_{ijij} +\overline{R}_{klkl} -2W_{ijij}+2H^2\notag\\
&+(\lambda_i-H)(\lambda_j-H) +(\lambda_k-H)(\lambda_l-H).\notag\\
=&\overline{R}_{ijij} +\overline{R}_{klkl} -2W_{ijij}+2H^2
+\mu_i\mu_j+\mu_k\mu_l.\notag
\end{align}
Note that the above formula has no summation on $i, j, k, l$.

In what follows, without lose of generality, we assume,
at the point $p$,
$$\mu_1\ge \mu_2\ge\mu_3\ge\mu_4.$$
Next we will split our proof of case (ii) into three parts
according to the values of $\mathcal{K}$
(defined in (\ref{lambda})) and $\mu_i$.

\begin{itemize}
\item[(a)] $\mathcal{K}(p)\ge 0$,
$\mu_1\ge \mu_2 \ge 0\ge\mu_3\ge\mu_4$. \\
Let $i=1$ and $j=2$ in (\ref{S}),
and since locally conformally flatness implies $W\equiv 0$,
we get
$$
\frac{S}{6} = \overline{R}_{ijij} +\overline{R}_{klkl} +2H^2
+\mu_1\mu_2
+\mu_3\mu_4 \ge 2H^2.
$$
Consequently, $S\ge 12H^2\ge -12(1+H^2).$
The remain proof of this part is similar with case (i).
\\
\item[(b)]
 $\mathcal{K}(p) < 0$ and
$\mu_1\ge \mu_2\ge\mu_3> 0> \mu_4$.\\
In this part, a direct computation gives
\begin{align}\label{mu3negative}
\sum_i\mu_i^3 = \sum_{i=1}^3 \mu_i - \left(\sum_{i=1}^3 \mu_i\right)^3
\le 0.
\end{align}
Without loss of generality, we assume $H(p)\ge 0$
(otherwise the term "$12H\sum_i\lambda_i^3$" in (\ref{Q1})
will be nonnegative, this case can be dealt with a similarly method
as next part (c)).
Therefore, combined (\ref{Q1}), (\ref{mu3negative}) and the assumption
$\mathcal{K}\le 0$, we obtain
\begin{align}\label{Q2}
Q\le 3\left(1+H^2\right)^2
 -\frac{1}{24}\vert \mathring{A}\vert^2\left( \vert \mathring{A}\vert^2
 - 12(1-H^2)
 \right).
\end{align}
Note that  in case (ii)
$$\vert \mathring{A}\vert^2
\ge 12 +24H^2 \ge 12(1-H^2).$$
Thus the second term in the right hand side of (\ref{Q2})
is nonpositive, and consequently we have $Q\le 3(1+H^2)^2$.
\\

\item[(c)] $\mathcal{K}(p) <0$ and $\mu_1 >0> \mu_2\ge\mu_3\ge \mu_4$.\\
In this part, inequality (\ref{Q1}) is not enough for our estimate,
and we will go back into the equality (\ref{SR}) and estimate term
by term.\\
Firstly,  $\sum_i\mu_i^3\ge 0$, and we will use the following
inequality (see \cite[Lemma 1]{XuXu})
\begin{align}\label{mu3-1}
\sum_i \mu^3 \le \frac{1}{\sqrt{3}}\vert \mathring{A}\vert^3.
\end{align}
Secondly, for the term $6\sum_i \mu_i^2\left(a_{ii}-\frac{\sigma}{3}\right)$,
under the assumption of this part, we will use
 a more accurate (than (\ref{R})) estimate as follows.
\begin{align}\label{R1}
&3\sum_i \mu_i^2\left(a_{ii}-\frac{\sigma}{3}\right)\\
=&\mu_1^2\left( \overline{R}_{1212} +\overline{R}_{1313}
+\overline{R}_{1414}-2\left(\overline{R}_{2323}+
\overline{R}_{2424}+\overline{R}_{3434}\right)\right)\notag\\
&+\mu_2^2\left( \overline{R}_{2121} +\overline{R}_{2323}
+\overline{R}_{2424}-2\left(\overline{R}_{1313}+
\overline{R}_{1414}+\overline{R}_{3434}\right)\right)\notag\\
&+\mu_3^2\left( \overline{R}_{3131} +\overline{R}_{3232}
+\overline{R}_{3434}-2\left(\overline{R}_{1212}+
\overline{R}_{1414}+\overline{R}_{2424}\right)\right)\notag\\
&+\mu_4^2\left( \overline{R}_{4141} +\overline{R}_{4242}
+\overline{R}_{4343}-2\left(\overline{R}_{1212}+
\overline{R}_{1313}+\overline{R}_{2323}\right)\right)\notag\\
=&(\mu_1^2+\mu_2^2-2(\mu_3^2+\mu_4^2))\overline{R}_{1212}
+(\mu_1^2+\mu_3^2-2(\mu_2^2+\mu_4^2))\overline{R}_{1313}\notag\\
&+(\mu_1^2+\mu_4^2-2(\mu_2^2+\mu_3^2))\overline{R}_{1414}
+(\mu_2^2+\mu_3^2-2(\mu_1^2+\mu_4^2))\overline{R}_{2323}\notag\\
&+(\mu_2^2+\mu_4^2-2(\mu_1^2+\mu_3^2))\overline{R}_{2424}
+(\mu_3^2+\mu_4^2-2(\mu_1^2+\mu_2^2))\overline{R}_{3434}\notag\\
\le& (\mu_1^2+\mu_2^2-2(\mu_3^2+\mu_4^2))\overline{R}_{1212}
+(\mu_1^2+\mu_3^2-2(\mu_2^2+\mu_4^2))\overline{R}_{1313}\notag\\
&+(\mu_1^2+\mu_4^2-2(\mu_2^2+\mu_3^2))\overline{R}_{1414}\notag\\
=&(-2\mu_1\mu_2-(\mu_3-\mu_4)^2)\overline{R}_{1212}
+(-2\mu_1\mu_3-(\mu_2-\mu_4)^2)\overline{R}_{1313}\notag\\
&+(-2\mu_1\mu_4-(\mu_2-\mu_3)^2)\overline{R}_{1414}\notag\\
\le & -2\mu_1(\mu_2+\mu_3+\mu_4)=2\mu_1^2
\le \frac{3}{2}\vert\mathring{A}\vert^2,\notag
\end{align}
where we have used the facts that
$$\mu_1> 0>\mu_2\ge\mu_3\ge\mu_4, \quad \sum_i\mu_i=0,$$
and the inequality
\begin{align*}
\vert\mathring{A}\vert^2=\sum_i\mu_i^2
 \ge \mu_1^2 +\frac{(\mu_2+\mu_3+\mu_4)^2}{3}=\frac{4}{3} \mu_1^2.
\end{align*}
Combined (\ref{SR}), (\ref{lambda}),  (\ref{HA}), (\ref{mu3-1}),
(\ref{R1}) and the fact that
$0\le\sigma\le 12$, we obtain
\begin{align}\label{Q3}
Q\le &3\left(1+H^2\right)^2\\
 &+\frac{1}{12}\left(-\frac{1}{2} \vert \mathring{A}\vert^4
+4\sqrt{3}\vert H\vert\vert\mathring{A}\vert^3 + 3(1-2H^2)\vert \mathring{A}\vert^2
  \right)\notag\\
  := &3\left(1+H^2\right)^2 +F(\vert\mathring{A}\vert),\notag
\end{align}
where
\begin{align}\label{F}
F(\vert\mathring{A}\vert) &= \frac{1}{12}\left(-\frac{1}{2} \vert \mathring{A}\vert^4
+4\sqrt{3}\vert H\vert \vert\mathring{A}\vert^3 + 3(1-2H^2)\vert \mathring{A}\vert^2
  \right)\\
  &= \frac{1}{12}\left(-\frac{1}{2}( \vert\mathring{A}\vert^4
  -6\vert\mathring{A}\vert^2)+\vert H\vert(4\sqrt{3}
  \vert\mathring{A}\vert^3 -6\vert H\vert \vert\mathring{A}\vert^2) \right).\notag
\end{align}
It is easy to see $F(x)$ attains its maximum at
$x_0 = 3\sqrt{3}\vert H\vert + \sqrt{3+21H^2}$ and decreasing
when $x\ge x_0$.
\\
Keep in mind that this part is one of the three parts of
case (ii), which assumes that
$$\vert\mathring{A}\vert^2 \ge
12+ 24H^2.$$
Therefore,
if $x_0\le \sqrt{12+24H^2}$,  we have
\begin{align}\label{F1}
F(\vert\mathring{A}\vert)&\le F(\sqrt{12+24H^2})\\
&= -12H^2 (2H^2+1) + \vert H\vert
\left.\left(\frac{\sqrt{3}}{3} x^3
-\frac{1}{2} \vert H\vert x^2\right)\right\vert_{x=\sqrt{12+24H^2}}\notag\\
 &\le \vert H\vert f_1(\vert H\vert),\notag
\end{align}
where
$$
f_1(\vert H\vert) =\left.\left(\frac{\sqrt{3}}{3} x^3
-\frac{1}{2} \vert H\vert x^2\right)\right\vert_{x=\sqrt{12+24H^2}}\ge 0.
$$
If $x_0\ge \sqrt{12+24H^2}$, we obtain
\begin{align}\label{F2}
F(\vert\mathring{A}\vert) &\le F(x_0)\\
&=\frac{1}{12}\left(-\frac{1}{2}( x^4-6x^2)\vert_{x=x_0}
+\vert H\vert(4\sqrt{3}x^3 -6\vert H\vert x^2)\vert_{x=x_0}
  \right)\notag
\\
& \le \frac{1}{12}\left(-\frac{1}{2}( x^4-6x^2)\vert_{x=\sqrt{12+24H^2}}\right)
+\vert H\vert f_2(\vert H\vert)
  \notag,\\
  &\le -12H^2 (2H^2+1) +\vert H\vert f_2(\vert H\vert)\notag\\
  &\le \vert H\vert f_2(\vert H\vert),\notag
\end{align}
where
$$
f_2(\vert H\vert) =\left.\left(\frac{\sqrt{3}}{3} x^3
-\frac{1}{2} \vert H\vert x^2\right)\right\vert_{x=x_0}\ge 0.
$$

Combined (\ref{Q3}), (\ref{F}), (\ref{F1}) and (\ref{F2}), we get
\begin{align}\label{Q4}
Q\le 3(1+H^2)^2 + 3\vert H\vert f(\vert H\vert),
\end{align}
where $f(\vert H\vert)$ is a function of $\vert H\vert$ defined by
\begin{equation}\label{f}
3f(\vert H\vert) := \left\{
\begin{aligned}
f_1(\vert H\vert), \quad x_0\le \sqrt{12+24H^2},\\
f_2(\vert H\vert), \quad x_0\ge \sqrt{12+24H^2}.
\end{aligned}
\right.
\end{equation}

\end{itemize}

To sum up the above two cases, we have proved the Claim, and
inequality (\ref{vol}) follows immediately.

Next we will show if $\vert H\vert$ is small, inequality
(\ref{vol1}) holds. Check all the cases in the proof of
the Claim, we find inequality (\ref{vol1}) holds except for the case
(ii) (c).  Thus, it is enough to show,
if $\vert H\vert$ is small, inequality (\ref{vol1}) holds
in the case (ii) (c). By (\ref{Q3}), it is sufficient to show
$F(\vert\mathring{A}\vert)\le 0$ when $\vert H\vert$ is small.
Observe that $F(\vert\mathring{A}\vert)$ can be decomposed as
\begin{align}\label{F3}
F(\vert\mathring{A}\vert)
=-\frac{\vert\mathring{A}\vert^2}{24}
\left(\vert\mathring{A}\vert -\eta_1\right)
\left(\vert\mathring{A}\vert -\eta_2\right),
\end{align}
where $\eta_1 = 4\sqrt{3}\vert H\vert -\sqrt{6+36H^2}$,
$\eta_2 = 4\sqrt{3}\vert H\vert+\sqrt{6+36H^2}$. Remeber that
in case (ii), $\vert\mathring{A}\vert\ge \sqrt{12+24H^2}$.
It is easy to see, if $\vert H\vert$ is small, say
$\vert H\vert \le\varepsilon_0$ for some constant
$\varepsilon_0$,
$$
\vert\mathring{A}\vert\ge \sqrt{12+24H^2} \ge \eta_2 >\eta_1,
$$
which implies
$F(\vert \mathring{A}\vert)\le 0$.

Check the above arguments step by step, we find the equality holds
in (\ref{vol1}) if and only
if
$$
\sigma = \sum_{i,j} \overline{R}_{ijij} = 12, \quad \mathring{A}\equiv 0,
$$
which implies $\Sigma$ is totally umbilic and
$$\overline{R}_{ijij} =1, \quad \text{for all}\quad i\neq j.$$
Therefore, by the Gauss equation, we get, for all $i\neq j$,
\begin{align}\label{Schur}
R_{ijij} = \overline{R}_{ijij} +\lambda_i\lambda_j
= 1 + (\mu_i+H)(\mu_j+H) =1+H^2,
\end{align}
which means the sectional curvature of $\Sigma$ at one point
$p$ is the same for all tagent plane $\pi\in T_p\Sigma$.
By Schur's lemma, $Sec_{\Sigma}$ is constant. Hence, by (\ref{Schur}),
$H$ is constant and $Sec_{\Sigma} \equiv 1+H^2$. Therefore, $\Sigma$ is
isometric to $\mathbb S^4(\frac{1}{1+H^2})$.\qed\\\\
{\bf Remark.} The condition "$\vert H\vert \le \varepsilon_0$"
is just a technical condition. The constant
$\varepsilon_0$ can be taken to be  $\sqrt{\frac{368\sqrt{3} - 598}{46}}$.
But this  is not the best number. Actually, after a long calculation
similar as (\ref{R1}), we can get a better estimate than (\ref{HA}) and
finally improve $\varepsilon_0$. We believe the condition
"$\vert H\vert \le \varepsilon_0$" is not necessary for inequality (\ref{vol1}).
Actually, in next section, we study a special case when $M$ is
rotationally symmetric and show that inequality (\ref{vol1}) holds for all $H$.

\section{A special case}
In this section, we will deal with a special case,
the ambient manifold is rotionally symmetric, i.e.,
$M= \mathbb R\times_\varphi \mathbb S^n$
with the metric
\begin{align}\label{metric}
g = dt^2 +\varphi^2(t) ds^2_n,
\end{align}
where $\varphi(t)$ is a smooth positive function, and
$ds^2_n$ is the standard metric of $\mathbb S^n$. Denote
by $\partial_t$ the unit vector in the $\mathbb R$ direction,
and assume
$X, Y$ are two vectors tangent to $\mathbb S^n$, then
the curvature tensor is given by (see \cite[section 4.2.3]{Petersen})
\begin{align}\label{peter}
\overline{R}(X\wedge \partial_t)
= -\frac{\ddot{\varphi}}{\varphi}X\wedge \partial_t, \quad
\overline{R}(X\wedge Y)
= \frac{1- \dot{\varphi}^2}{\varphi^2}X\wedge Y.
\end{align}
For simplicity we write
\begin{align}\label{kappa}
\kappa_1 := -\frac{\ddot{\varphi}}{\varphi}, \quad
\kappa_2 := \frac{1- \dot{\varphi}^2}{\varphi^2}.
\end{align}
Let $\Sigma$ be a hypersurfaces in $M$,
$T$ be the tangential (with $\Sigma$) part of $\partial_t$,
$e_1, \cdots, e_n$ be the local orthonormal frame on $\Sigma$.
Write $T_i = g(T, e_i)$.
Decomposed each $e_i$ into two parts
\begin{align}\label{decom}
e_i = e_i^\prime + g( e_i, \partial_t)
= e_i^\prime + T_i,
\end{align}
where $e_i^{\prime}$ tantgent to $\mathbb S^n$.
A direct computation, by using (\ref{peter}), (\ref{decom}) and the multilinerity
of the curvature tensor, gives
\begin{align}\label{gauss}
\overline{R}_{ijkl}=
\kappa_2(\delta_{ik}\delta_{jl} - \delta_{il}\delta_{jk})
+(\kappa_1-\kappa_2)\left( T_iT_k\delta_{jl}+T_jT_l\delta_{ik}
-T_iT_l\delta_{jk} - T_jT_k\delta_{il} \right).
\end{align}
Therefore, for $i\neq j$, we have
\begin{align}\label{gauss1}
R_{ijij}=\overline{R}_{ijij}+\lambda_i\lambda_j=
\kappa_2
+(\kappa_1-\kappa_2)\left( T_i^2+T_j^2 \right)+\lambda_i\lambda_j.
\end{align}

We need the following lemma, which was first proved by Cartan (we appreciate
Professor Marcos Dajczer pointing this fact out to us). For completeness,
we give a direct proof here.
\begin{lemma}\label{principal}
Let $M$ be an $n+1 (n\ge 4)$ dimensional rotationally symmetric
Riemannian manifold with metric (\ref{metric}), and $\Sigma$ be a
hypersurface in $M$.  Then $\Sigma$ is locally conformally flat
if and only if at each point $p\in \Sigma$, there
are at most two distinct principal curvatures, one of them has multiplicity
$n-1$.
\end{lemma}
\noindent
{\it Proof.}
We will adopt the notations in section 2.
In this proof, $i, j, k, l$ will range from $1$ to $n$.
By using the Weyl
tensor formula (\cite[p117]{Aubin}), we have
\begin{align}\label{W}
W_{ijij} = &R_{ijij} - \frac{1}{n-2}
\left(R_{ii}
+ R_{jj}\right)+ \frac{S}{(n-1)(n-2)}\\
=&\kappa_2
+(\kappa_1-\kappa_2)\left( T_i^2+T_j^2 \right)
+\lambda_i\lambda_j\notag\\ &
-\frac{1}{n-2}
\sum_{k\neq i} \left(\kappa_2
+(\kappa_1-\kappa_2)\left( T_i^2+T_k^2 \right) +\lambda_i\lambda_k\right)\notag\\
&-\frac{1}{n-2} \sum_{k\neq j} \left(\kappa_2
+(\kappa_1-\kappa_2)\left( T_j^2+T_k^2 \right) +\lambda_j\lambda_k\right)\notag\\
&+ \frac{\sum_{k\neq l} \left(\kappa_2
+(\kappa_1-\kappa_2)\left( T_k^2+T_l^2 \right) \right)+n^2H^2 -\vert A\vert^2}{(n-1)(n-2)}\notag\\
=&\lambda_i\lambda_j -\frac{1}{n-2}
\left(\sum_{k\neq i}  \lambda_i\lambda_k
+ \sum_{k\neq j} \lambda_j\lambda_k\right)
+ \frac{n^2H^2 -\vert A\vert^2}{(n-1)(n-2)}.\notag
 \end{align}
Using the relations
$$\lambda_k = \mu_k +H,\quad \vert A\vert^2
=\vert \mathring{A}\vert^2+nH^2,
$$
we
 substitute $\mu_k$ for $\lambda_k$ in the above
 equality and obtain
 \begin{align*}
W_{ijij} = \frac{(\mu_i+\mu_j)^2 +(n-4)\mu_i\mu_j}{n-2}
-\frac{\vert \mathring{A}\vert^2}{(n-1)(n-2)}.
 \end{align*}
Therefore, $\Sigma$ is
\begin{align*}
\text{locally conformally flat} &\Longleftrightarrow
W\equiv 0\\
(\text{at each point})&\Longleftrightarrow
(\mu_i+\mu_j)^2 +(n-4)\mu_i\mu_j
=\frac{\vert \mathring{A}\vert^2}{(n-1)}, \forall i\neq j\\
(\sum_i\mu_i = 0, \vert\mathring{A}\vert^2
= \sum_i\mu_i^2)&\Longleftrightarrow
\mu_i = \mu_j\ \ \text{or}\ \  \mu_i = -(n-1)\mu_j,
\forall i\neq j.
\end{align*}
Thus $\{\mu_1, \cdots, \mu_n\} = \{\mu, -(n-1)\mu\}$, and $\mu$
has multiplicity $n-1$. Consequently,
$$
\{\lambda_1, \cdots, \lambda_n\} = \{\mu+H, -(n-1)\mu+H\},
$$
and $\mu+H$ has multiplicity $n-1$.
\qed
%\\\\
%{\bf Remark.} Note that given special $\varphi$ in (\ref{metric}),
%$M$ will be of constant sectional curvature. Actually,
%\begin{equation*}
%M=\left\{
%  \begin{aligned}
%    &\mathbb S^{n+1}, \quad \varphi(t) = \sin(t),\\
%    &\mathbb R^{n+1}, \quad \varphi(t) = t,\\
%    &\mathbb H^{n+1}, \quad \varphi(t) = \sinh(t).
%  \end{aligned}
%  \right.
%\end{equation*}
%Therefore, the above lemma is a generalization of Cartan's
%result (see \cite{Cartan})
%which states that a hypersurface in a space form (dim$>4$)
%is locally conformally flat if and only if
%$n-1$ principal curvatures are the same.

With the aid of the above lemma, we can prove the following theorem.

\begin{theorem}\label{rotsym}
Let $M$ be a rotationally symmetric Riemannian 5-manifold with
$0\le \kappa_1 \le \kappa_2 \le 1$ ($\kappa_1$ and $\kappa_2$ are defined
in (\ref{kappa})),
and $\Sigma$ be a locally conformally flat closed hypersurface embedded
in $M$ with mean curvature
 $H$, then
$$
\int_{\Sigma} (1+H^2)^2 \ge \frac{4\pi^2}{3}\chi(\Sigma).
$$
The equality holds if and only if  $H$ is
constant, $\Sigma$ is
totally umbilic and isometric to $\mathbb S^4\left(\frac{1}{1+H^2}\right)$.
\end{theorem}
\noindent
{\it Proof.} We will adopt the same notations as in the proof of Theorem \ref{1}.
By Lemma \ref{principal}, we have
$\vert \mu_1\vert = 3\vert \mu_i\vert$, $i=2,3,4$. Therefore,
$$
\sum_i \mu_i^4 = \frac{7}{12}\vert \mathring{A}\vert^4.
$$
Direct computations by using  (\ref{gauss}) yield
\begin{align*}
  \sigma &= 12\kappa_2+6(\kappa_1-\kappa_2)\vert T\vert^2,\\
  a_{ii}&=3\kappa_2+(\kappa_1-\kappa_2)(2T_i^2+\vert T\vert^2),\\
  \vert\mathring{a}\vert^2&= 3(\kappa_1 - \kappa_2)^2\vert T\vert^4.
\end{align*}

Insert the above equalities into (\ref{SR}),
we obtain,
\begin{align}\label{Q5}
Q
 =&3(\kappa_2+H^2)^2+\frac{1}{12}\left( -\frac{3}{4}\vert \mathring{A}\vert^4
  +12H\sum_i \mu_i^3 -6(\kappa_2+3H^2)\vert \mathring{A}\vert^2\right)\\
  &+\frac{(\kappa_1-\kappa_2)}{2} \left(\left(6+\vert\mathring{A}\vert^2
  +4H^2\right)\vert T\vert^2
  +2\sum_i T_i^2(\mu_i-H)^2\right)\notag\\
 \le &3(\kappa_2+H^2)^2+\frac{1}{12}\left( -\frac{3}{4}\vert \mathring{A}\vert^4
  +12H\sum_i \mu_i^3 -6(\kappa_2+3H^2)\vert \mathring{A}\vert^2\right)\notag\\
  \le &3(\kappa_2+H^2)^2+\frac{1}{12}\left( -\frac{3}{4}\vert \mathring{A}\vert^4
  +4\sqrt{3}\vert H\vert\ \vert\mathring{A}\vert^3
  -6(\kappa_2+3H^2)\vert \mathring{A}\vert^2\right)\notag\\
  \le &3(\kappa_2+H^2)^2\le 3(1+H^2)^2.\notag
\end{align}
The remain proof is similar as the proof of Theorem \ref{1}.\qed
\\\\
{\bf Remark.} The assumption "$0\le \kappa_1\le\kappa_2\le 1$" is reasonable for
many manifolds. For example, \\
$\bullet$ if we take $\varphi(t) = \sin(t)$, then $M=\mathbb S^5$ and
$\kappa_1 = \kappa_2 \equiv 1$;\\
$\bullet$ if we take $\varphi(t)\equiv 1$, then $M = \mathbb S^4 \times \mathbb R$
and $0\equiv \kappa_1 < \kappa_2 \equiv 1.$

\end{document}